\documentclass[11pt]{amsart}
\usepackage{lipsum}
\usepackage{amsmath, amsthm, amsfonts, ifpdf}
\usepackage[dvipsnames,usenames]{color}
\usepackage{graphics,color}
\usepackage{multicol}
\usepackage{tikz}
\usepackage{pgfmath}
\usetikzlibrary{arrows,shapes,trees,backgrounds}
\usetikzlibrary{intersections}

\theoremstyle{plain}
\newtheorem{theorem}{Theorem}[section]
\newtheorem*{theorem*}{Theorem}

\newtheorem{pro}[theorem]{Proposition}
\newtheorem{Def}[theorem]{Definition}
\newtheorem{lem}[theorem]{Lemma}
\newtheorem{cor}[theorem]{Corollary}
\newenvironment{customthm}[1]
  {\renewcommand\theinnercustomthm{#1}\innercustomthm}
  {\endinnercustomthm}
  \newtheorem{innercustomthm}{Theorem}

\theoremstyle{definition}
\newtheorem*{Def*}{Definition}

\theoremstyle{remark}

\newtheorem{Rem}[theorem]{Remark}

\numberwithin{equation}{section}

\newcommand{\bpo}{\begin{pro}}
\newcommand{\epo}{\end{pro}}
\newcommand{\be}{\begin{equation}}
\newcommand{\ene}{\end{equation}}
\newcommand{\br}{\begin{Rem}}
\newcommand{\er}{\end{Rem}}
\newcommand{\bl}{\begin{lem}}
\newcommand{\el}{\end{lem}}
\newcommand{\bd}{\begin{Def}}
\newcommand{\ed}{\end{Def}}
\newcommand{\ben}{\begin{enumerate}}
\newcommand{\een}{\end{enumerate}}
\newcommand{\bp}{\begin{proof}}
\newcommand{\ep}{\end{proof}}
\newcommand{\beq}{\begin{equation*}}
\newcommand{\eeq}{\end{equation*}}
\newcommand{\bear}{\begin{eqnarray*}}
\newcommand{\eear}{\end{eqnarray*}}
\newcommand{\bt}{\begin{theorem}}
\newcommand{\et}{\end{theorem}}
\newcommand{\bst}{\begin{split}}
\newcommand{\est}{\end{split}}

\newcommand{\bal}{\begin{aligned}}
\newcommand{\eal}{\end{aligned}}
\newcommand{\bas}{\begin{align*}}
\newcommand{\eas}{\end{align*}}

\renewcommand{\P}{\partial}
\newcommand{\F}[2]{\frac{#1}{#2}}
\newcommand{\la}{\langle}
\newcommand{\ra}{\rangle}

\newcommand{\R}{\mathbb{R}}

\newcommand{\bnb}{\bar{\nabla}}
\newcommand{\nb}{\nabla}

\newcommand{\vn}{\P_r}

\newcommand{\Sc}{\varepsilon}

\newcommand{\vp}{\varphi}

\newcommand{\gu}{graph$(u)$}

\newcommand{\PLH}{{\mkern-1mu\times\mkern-1mu}}
\newcommand{\WP}{N\PLH_{\phi}\R}

\begin{document}
\title[mean curvature type flow of closed graphs]{Mean curvature type flows of graphs in product manifolds}
\author{Aijin Lin, Hengyu Zhou}
\address [Aijin Lin]{College of Science, National University of Defense Technology, Changsha 410073, P. R. China}
\email{aijinlin@pku.edu.cn}
\address[Hengyu Zhou]{Department of Mathematics, Sun Yat-sen University, 510275,
Guangzhou, P. R. China}
\email{zhouhy28@mail.sysu.edu.cn}

\date{\today}
\begin{abstract}
	 \indent In this note we study a large class of mean curvature type flows of graphs in product manifold $N\PLH\R$ where $N$ is a closed Riemannian manifold. Their speeds are the mean curvature of graphs plus a prescribed  function. We establish long time existence and uniformly convergence of those flows with a barrier condition and a condition on the derivative of prescribed function with respect to the height. As an application we construct a weighted mean curvature flow in large classes of warped product manifolds which evolves each graph into a totally geodesic slice. 
\end{abstract}
\maketitle
 \section{Introduction}
 In this paper we are interested in mean curvature type flows of closed graphs in product manifolds. \\
 \indent The first motivation of this project is to find existence of hypersurfaces with prescribed mean curvature in a Riemannian manifold. Suppose $M$ is a Riemannian manifold and $f(x,X)$ is a smooth function on $M\PLH TM$ where $TM$ denotes the tangent bundle of $M$. It is natural to ask whether we can find a smooth hypersurface $\Sigma$ satisfying 
    \be\label{basic:equation}
    H(x)=f(x,\vec{v})\quad x\in \Sigma
    \ene 
    where $\vec{v}$ is the normal vector of $\Sigma$
 (see Yau \cite{Yau82}). The solution to problem \eqref{basic:equation} can be viewed as a stationary point of the following nonlinear mean curvature type flow:
 \be \label{basic:parabolic}
\P_t F=-(H-f)\vec{v}
\ene 
where $H$ is the mean curvature of smooth hypersurfaces $F(t, .)$ with some initial condition.\\
\indent In general we need some geometric conditions to guarantee the existence of problem \eqref{basic:equation} and the long time existence of the flow in \eqref{basic:parabolic}. For the former one we give a concept referred as a barrier condition. For the latter one we will work with the flows of graphs in product manifolds.\\
\indent Suppose $\Omega$ is a smooth domain such that its boundary $\P\Omega$ has two connected components $\P_1\Omega$ and $\P_2\Omega$. Let $H_{\P_1\Omega}$ ($H_{\P_2\Omega})$ be the mean curvature of $\P_1\Omega$ ($\P_2\Omega$) with respect to the outward (inward) normal vector $\vec{v}_1$ ($\vec{v}_2$). The choice of the normal vector guarantees that the mean curvature of the boundary of a ball in Euclidean space is always a positive constant. A function $f(x,X)\in C^1(\bar{\Omega}\times T\bar{\Omega})$ is said to satisfies the barrier condition on $\Omega$ provided that 
\be\label{condition:label}
   H_{\P\Omega_1}\geq f(x,\vec{v}_1)\quad H_{\P\Omega_2} \leq  f(x,\vec{v}_2)
\ene 
\indent The barrier condition \eqref{condition:label} plays an essential role for the exsitence of hypersurfaces with prescribed mean curvature. By assuming this condition, in Euclidean space Treibergs-Wei \cite{TW83} obtained the existence of starshaped hypersurfaces with prescribed mean curvature function when the function  $f(p)$ satisfies a barrier condition and $\P_u(uf(x,u))\leq 0$ where $x\in S^n$ and $p=(x,u)$. Zhu \cite{Zhu94} established a similar existence result in Euclidean space via a mean curvature type flow when $f(p)$ is a concave function and satisfies a barrier condition. de Andra-Barbosa-de
 Lira  \cite{ABL09} and Cheng-Li-Wang \cite{CLW17} considered the existence of hypersurfaces with general curvature function in warped product manifolds. Their conditions include a barrier condition and $\P_u(\phi(u)f(x,u))\leq 0$ where $\phi$ is a warped function and $(x,u)$ denotes the coordinates of warped product manifolds (see Section 2.3). \\
 \indent Another motivation of this project is the nonparametric mean curvature type flow with prescribed contact angle which is a nonparametric version of \eqref{basic:parabolic}. Such flows with prescribed contact angle were studied by Huisken \cite{Hui89A} in the case of mean curvature flow and  by  Guan \cite{GB96} in the case of general  nonparametric mean curvature type flows of graphs over a bounded domain in Euclidean space. A central idea in \cite{GB96} is that the $C^0$ estimate will imply the $C^1$ estimate in mean curvature type equation \eqref{basic:equation} provided that $f$ satisfies an admissible condition. This method was recently developed in general Riemannian manifolds by the second author in \cite{Zhou17}. \\
    \indent In this note we continue to apply the idea of Guan \cite{GB96} into the nonparametric mean curvature type flow without boundary. The main result of this paper is stated as follows. 
    \bt\label{thm:main}  Suppose $N$ is a closed manifold and  $h(x,u), g(x,u)$ be smooth functions on product manifold $N\PLH\R$. Suppose 
    	\begin{enumerate}
    		\item For all $x\in N$ there exists two constants $u_0<u_1$ satisfying 
    		$$
    		g(x,u_0)+h(x,u_0)\geq 0\quad  g(x,u_1)+h(x,u_1)\leq 0
    	     $$
    		\item Moreover $\P_u g(x,u)\leq 0$ for all $x\in N$ and $u\in [u_0,u_1]$. 
    	\end{enumerate} 
    	If $\Sigma_0$ is a smooth graph contained in the domain $N\PLH [u_0, u^0]$, then the flow 
    	\be \label{mct:flow}
    	\P_t F=-(H-h(x,u)\la \vec{v},\P_r\ra +g(x,u))\vec{v}
    	\ene
    	with initial smooth graph $\Sigma_0$ in $N\PLH\R$ exists for all time, remains graphical. If $h(x,u)=h(u)$, $u_t(x,t)$ converges to $0$ uniformly.  Moreover there is a $\{t_i\}_{i=1}^\infty$ such that $t_i\rightarrow \infty$ and $u(x,t_i)$ converges uniformly to a smooth solution to problem $$ H(u)=h(u)\la \vec{v},\P_r\ra  -g(x,u)$$
    	where $H(u)$ is the mean curvature of the graph of $u(x)$ in $N\PLH\R$. In the above $\vec{v}$ always denotes the downward normal vector. 
    \et
     \br  For mean curvature flows in various warped product manifolds, we refer to the work of Borisenko-Miquel \cite{BM12} and Huang-Zhang-Zhou \cite{Zhou17}. \\
     \indent Condition (1) is a special case of the barrier condition \eqref{condition:label} and will give the $C^0$ estimate of the flow \eqref{mct:flow} (see Lemma \ref{stepone}). Thus we can apply the idea of Guan \cite{GB96}.\\
     \indent If choose general $h(x,u)$ depending on $x$, we can not obtain the integral estimate in Lemma \ref{lm:integral}. This is a key estimate to show that $u_t$ converges to $0$ uniformly in Theorem \ref{thm:main}. 
        \er 
   \indent According to Lemma \ref{lm:transform} there is an one-to-one correspondence between the graphs in product manifolds and those in warped product manifolds (see Section 2.3). Therefore the uniform convergence in Theorem \ref{thm:main} gives a flow-type proof of the following existence result. 
    \begin{cor} \label{cor:one} Suppose $N$ is a n-dimensional closed Riemannian manifold. Let $\phi(u)$ be a positive smooth function on $\R$ and $f(x,u)$ be a smooth function where $x\in N$ and $u\in \R$. Suppose the following two conditions hold:
       \begin{enumerate}
        	\item  For all $x\in N$ there exists two constants $u_0<u_1$ satisfying 
        \be\label{con:bar}
        f(x,u_0)\geq n\F{\phi'(u_0)}{\phi(u_0)}\quad  f(x,u_1)\leq n\F{\phi'(u_1)}{\phi(u_1)}
        \ene
        \item Moreover $\P_u (f(x,u)\phi(u))\leq 0$ for all $x\in N$ and $u\in [u_0,u_1]$. 
        \end{enumerate}
     Then there exists a  smooth hypersurface $\Sigma=(x,u(x))$ in warped product manifold $N\PLH_{\phi}\R$ such that its mean curvature is $f(x,u)$. 
    \end{cor}
\br \eqref{con:bar} is the special form of the  barrier condition \eqref{condition:label} in warped product manifolds. \\
\indent  Here the sign of $\phi'(u)$ can change comparing to the existence results in \cite{CLW17} and \cite{ABL09}. This fact is useful when we apply this result into Fuchsian manifolds where $\phi(u)=\cosh(u)$ for $u\in \R$ and $N$ is a closed hyperbolic surface (see \cite{HW13}). 
\er

 At last we find some good behaviors of a weighted mean curvature flow in warped product manifolds according to Theorem \ref{thm:main}. It may have independent interests. 
 \begin{cor}\label{cor:two} Let $a<b$ be two fixed constants. Suppose a warped function $\phi(u)$ satisfies
 	\be\label{condition:a}
 	  \phi'(a)\leq 0\quad \phi'(b)>0 \quad \phi''(u)\geq 0\quad \text{for all $u\in (a,b)$}
 	\ene
 	Let $N$ be a closed Riemananian manifold. Suppose $\Sigma_0$ is a graphical smooth hypersurface in a warped product manifold $N\PLH_{\phi}(a,b)$. Consider a weighted mean curvature flow 
 	 \be \label{wt:flow}
 	  \P_t F(p)=-\phi^2(u)H\vec{v}
 	 \ene 
 	 where $p=(x,u)\in N\PLH_{\phi}\R$, $F(0,p)=p$ for $p\in \Sigma$, $H$ is the mean curvature of $F_t(\Sigma_0)$. Then the flow of $\Sigma_0$ in \eqref{wt:flow} exists for all times, remains graphical and converges uniformly to the totally geodesic slice in $N\PLH_{\phi}[a,b]$. 
 	\end{cor}
 \br The mean curvature flow  in warped product manifolds can blow up in finite time see Appendix A of \cite{Zhou17}. In this note such flow is referred as a weighted mean curvature flow. \\
 \indent There are many important warped product manifolds satisfying \eqref{condition:a}. We refer the readers to \cite{Bre13}. 
 \er 
 \indent  This paper is organized as follows.  In Section 2 we discuss some preliminary facts about product manifolds and warped product manifolds. We also collect necessary evolution equations on nonparametric mean curvature flow in product manifolds from Section 3 in \cite{Zhou17}. In Section 3 we establish Theorem \ref{thm:main}. In Section 4 we demonstrate Corollary \ref{cor:one} and Corollary \ref{cor:two}. 
\section{Preliminary facts}
In this section we collect some facts on product manifolds, warped product manifolds and nonparametric mean curvature type flows. For more details we refer to \cite{Zhou17} and \cite{Zhou16a}.
\subsection{The graphs in product manifolds}Throughout this paper $N$ is a $n$-dimensional closed Riemannian manifold with Riemannian metric $\sigma$. Let $u(x)$ be a smooth function over $N$. Let $\{\P_i:=\F{\P}{\P x_i}\}$ be a local frame along $N$ and $\{dx^i\}$ be its dual frame. Then the metric $\sigma$ takes the form 
\be 
\sigma=\sigma_{ij}dx^idx^j
\ene
with $(\sigma^{ij})=(\sigma_{ij})^{-1}$. Let $u_{i}$ and $u_{jk}$ denote the covariant derivatives with respect to $\sigma$. The gradient of $u$ is 
 \be 
 Du=u^k\P_k \quad u^k=\sigma^{kl}\P_l
 \ene 
Let $\Sigma$ denote the graph of $u(x)$ in the product manifold $N\PLH\R$. Then its downward normal vector is
\be \label{vec:graph}
\vec{v}=\F{Du-\P_r}{\omega}
\ene 
where $\omega =\sqrt{1+|Du|^2}$.  Let $X_j$ denote $\P_j+u_j\P_r$. Then induced the metric of $\Sigma$ is 
 $$
 g_{ij}=\sigma_{ij}+u_ju_i
 $$
 with the inverse matrix $(g^{ij})=(g_{ij})^{-1}$ 
 $$g^{ij}=\sigma^{ij}-\F{u^iu^j}{\omega^2}$$
 The second fundamental form of {\gu} is
 \be
 h_{ij}=\la \bnb_{X_i}\vec{v}, X_j\ra =\F{u_{kj}}{\omega};
 \ene
 where $\bnb$ is the covariant derivative of $N\PLH\R$. The mean curvature of $\Sigma$ is 
 \be\label{mc:graph}
  H=\bar{d}iv(\vec{v})=g^{ij}h_{ij}=\F{1}{\omega}g^{ij}u_{ij}
 \ene 
 where $\bar{d}iv$ is the divergence of $N\PLH\R$. The norm of the second fundamental form is
 \be
 |A|^2=g^{kl}g^{ij}h_{ki}h_{lj} =\F{g^{il}g^{kj}u_{kl}u_{ij}}{\omega^2};
 \ene
\subsection{Nonparametric form}Now we consider the nonparametric form of the mean curvature type flow in  \eqref{mct:flow}. 
\bl\label{lm:nonform} Suppose $\{\Sigma_t\}$ is a family of graphs over $N$ satisfying \eqref{mct:flow} and is represented as $\{(x,u(x,t))\}$ where $u(x,t):N\times [0,T]\rightarrow \R$. Then $u(x,t)$ satisfies the following quasilinear parabolic equation:
 \be \label{equ:non:form}
 \P_t u=g^{ij}u_{ij}+h(x,u)+g(x,u)\omega
 \ene 
\el 
\bp According to \eqref{vec:graph} $\la\vec{v},\P_r\ra=-\F{1}{\omega}$. Then according to \eqref{mc:graph} $F'(x,t)=(x,u(x,t))$ satisfies that 
\be 
   (\P_t F'(x,t))^{\bot}=-\F{1}{\omega}(g^{ij}u_{ij}+h(x,u)+g(x,u)\omega)
\ene
where $\bot$ is the projection into the normal bundle of $\Sigma_t$
 By definition 
   $$
    (\P_t F'(x,t))^{\bot}=\P_t u\la \vec{v},\P_r\ra =-\P_t u\F{1}{\omega}
   $$
   Therefore we obtain the lemma. 
\ep
Now we define 
  $$
  \psi(x,u,Du)=-h(x,u)-g(x,u)\omega
  $$
  with $\psi_k :=\F{\P\psi(x,u,Du)}{\P x_k}$ and $\psi_{u_k}:=\F{\P \psi}{\P u_k}(x,u,Du)$. Here we view $Du$ as a function of $u^k$ written as $Du=u^k\P_k$. \\
 \indent Therefore we proceed as in Section 3 of \cite{Zhou17} and define a parabolic operator 
  \be \label{def:operator}
  L=g^{ij}\nb_i\nb_j-\psi_{u_k}\nb_k-\P_t
\ene   
Now we record the evolution equations of $u(x,t)$ and $\omega=\sqrt{1+|Du|^2}$ along the flow \eqref{equ:non:form} as follows. 
\bl Let $u(x,t)$ be the smooth solution to the flow \eqref{equ:non:form} on $[0,T]$.  Then $\omega$ satisfies that 
\begin{gather}
Lu =-\psi-\psi_{u_k}\nb_k u\\
L\omega=(|A|^2+Ric(v_N,v_N)+\psi_u\la v_N,v_N\ra)\omega+\F{2}{\omega}g^{il}\omega_i\omega_l +\F{u^k}{\omega}\psi_k \label{met}
\end{gather}
where $\psi=-h(x,u)-g(x,u)\omega$ and $v_N=\F{Du}{\omega}$. Here $Ric$ denotes the Ricci curvature of $N$. 
\el 
\bp The first identity is obvious according to \eqref{equ:non:form}. For the proof of the second identity see Lemma 3.5 in \cite{Zhou17}. 
\ep 
\subsection{Warped product manifolds} Let $I$ be an open interval and $\phi(r)$ be a smooth positive function on $I$. Recall that a warped product manifold $N\PLH_{\phi}I$ is defined the set $\{(x,r):x\in N, r\in I\}$ equipped with the metric $\phi^2(r)\sigma+dr^2$.\\
\indent  The warped product manifold is a generalization of space forms. When $N$ is the standard sphere and $I$ is $(0,\infty)$, these warped product manifolds are called as rotationally symmetric spaces (see Scheuer \cite{Sch17}). When $N$ is a closed hyperbolic surface and $\phi(r)=\cosh(r)$, such warped product manifold is referred as a Fuchsian manifold (see Huang-Wang \cite{HW13}). For more examples of warped product manifolds, we refer to Brendle \cite{Bre13}.  \\
\indent  Now we study graphs in warped product manifolds and establish their connection with those in product manifolds. \\
\indent Let $u(x)$ be a smooth function on $N$ and $\Sigma$ be the graph of $u(x)$ in $\WP$. The downward normal vector to $\Sigma$ is given by
\be\label{eq:normal}
\vec{v}=\F{Du-\phi^2(r)\P_r}{\sqrt{\phi^2(r)|Du|^2+\phi^4(r)}}
\ene
where $|Du|^2$ is the norm of the gradient $Du$ in $(M,\sigma)$ and $r=u(x)$. We define
\be\label{relation}
\vp(x)=\Phi(u(x)), \quad u(x)=\Phi^{-1}(\vp(x))
\ene
where $\Phi(r)$ is a strictly increasing function with $\Phi'(r)=\F{1}{\phi(r)}$ and $\Phi^{-1}(r)$ is its inverse function. Here $\Phi(r)$ may be negative. With these notation the downward normal vector in \eqref{eq:normal} is rewritten as
\be
\vec{v}=\F{\F{D\vp}{\phi(r)}-\P_r}{\omega}
\ene
where $\omega=\sqrt{1+|D\vp|^2}$ and $r=\Phi^{-1}(\vp(x))$. The mean curvature vector on $\Sigma$ is $-H\vec{v}$ where $H$ is the mean curvature of $\Sigma$. Let $\sigma=\sigma_{ij}dx_i dx_j$. Then a direct computation shows that
\bl\label{lm:ef} The mean curvature of $\Sigma=(x,u(x))$ in $\WP$ is
\be
H=\F{1}{\omega\phi(r)}((\sigma^{ij}-\F{\vp^i\vp^j}{\omega^2})\vp_{ij}-n\phi'(r))
\ene
where $r=u(x)$, $\vp(x)=\Phi(u(x))$, $(\sigma^{ij})=(\sigma_{ij})^{-1}$ and $n$ is the dimension of $M$.
\el
\bp Let $\bnb$ and $\la\ \ra$ denote the covariant derivative and the inner product of $N\PLH_{\phi}\R$ respectively. Let $\{\P_i\}_{i=1}^n$, $\P_r$ be a local frame of $N$ and $\R$ respectively. Moreover, $\sigma=\sigma_{ij}dx^idx^j$. \\
\indent All computations are evaluated at $(x,r)$. A well-known fact in warped product manifolds is that $V=\phi(r)\P_r$ is a conformal vector in $\WP$ satisfying $\bnb_{X}V=\phi'(r)X$. In particular
\be
\bnb_{\P_i}\P_r=\bnb_{\P_r}\P_i=\F{\phi'(r)}{\phi(r)}\P_i
\ene
On the other hand $\{X_i:\P_i+ u_i(x)\P_r\}$ is a local frame on $\Sigma$. The first fundamental form of $\Sigma$ is $$
g_{ij}=\la X_i, X_j\ra=\phi^{2}(r)\{\sigma_{ij}+\vp_i\vp_j\}
$$
Thus its inverse is $g^{ij}=\phi^{-2}(r)(\sigma^{ij}-\vp^i\vp^j\omega^{-2})$ where $\omega=\sqrt{1+|D\vp|^2}$. The second form follows from
\begin{align*}
h_{ij}&=\la\bnb_{(\P_i+\vp_i V)}\vec{v},\P_j+\vp_j V\ra;\\
&=-\F{\phi(r)}{\omega}(\vp_{ij}-\phi'(r)(\sigma_{ij}+\vp_i\vp_j));
\end{align*}
In the last line above some hidden facts are
$$\la\bnb_{\P_i}\P_j, \P_k\ra=\phi^2(r) \la\nb_{\P_i}\P_j, \P_k\ra_\sigma\quad \la\bnb_{\P_i}\P_j, \vn\ra=\phi'(r)\phi(r)\sigma_{ij}$$ where $\nb$ and $\la \, \ra_{\sigma}$ are the covariant derivative and the inner product of $(M,\sigma)$ respectively. The conclusion follows from $H=g^{ij}h_{ij}$.
\ep
A direct application of Lemma \ref{lm:ef} establishes a correspondence between graphs in warped product manifolds and those in product manifolds. A similar but slightly different case is discussed in Gerhardt \cite{GH98}.
\bl\label{lm:transform} Let $u(x)$ be a smooth function on $N$. The mean curvature of the graph of $u(x)$ in $N\PLH_\phi\R$ is $f(x,u(x))$ if and only if
\be\label{eq:rst}
H(\vp)=f(x,u)\phi(u)+n\F{\phi'(u)}{\omega}
\ene
where $u(x)=\Phi^{-1}(\vp(x))$ and $\omega=\sqrt{1+|D\vp|^2}$. Here $H(\vp)$ is the mean curvature of the graph $\vp(x)$ in product manifold $N\PLH\R$. 
\el
\bp The proof is straightforward by combining Lemma \ref{lm:ef} with \eqref{mc:graph}. 
\ep

\section{The proof of Theorem \ref{thm:main} }
In this section we give the proof of Theorem \ref{thm:main}. Comparing to  the prescribed contact angle conditions in \cite{Zhou17} here we discuss mean curvature type flows without boundary condition in product manifolds. \\
\indent For the convenience of readers Theorem \ref{thm:main} is restated as follows. 
    \begin{customthm}{\ref{thm:main}} Suppose $N$ is a closed manifold and  $h(x,u), g(x,u)$ be smooth functions on $N\PLH\R$. Suppose 
 \begin{enumerate}
		\item For all $x\in N$ there exists two constants $u_0<u_1$ satisfying 
		$$ 
		g(x,u_0)+h(x,u_0)\geq 0\quad  g(x,u_1)+h(x,u_1)\leq 0
		$$
		\item Moreover $\P_u (g(x,u))\leq 0$ for all $x\in N$ and $u\in (u_0,u_1)$. 
	\end{enumerate} 
 If $\Sigma_0$ is a smooth graph contained in the domain $N\PLH (u_0, u_1)$, then the flow 
$$
\P_t F=-(H-h(x,u)\la \vec{v},\P_r\ra +g(x,u))\vec{v}
$$
with initial smooth graph $\Sigma_0$ in $N\PLH\R$ exists for all time, remains graphical.
 If $h(x,u)=h(u)$, $u_t(x,t)$ converges to $0$ uniformly.  Moreover there is a $\{t_i\}_{i=1}^\infty$ such that $t_i\rightarrow \infty$ and $u(x,t_i)$ converges uniformly to a smooth solution to problem $$ H(u)=h(u)\la \vec{v},\P_r\ra  -g(x,u)$$
where $H(u)$ is the mean curvature of the graph of $u(x)$ in $N\PLH\R$. In the above $\vec{v}$ always denotes the downward normal vector. 
\end{customthm}    
 \indent Our proof is divided into three steps: the $C^0$ estimate, the $C^1$ estimate and the uniformly convergence. In subsection 3.1 we establish the long time existence in Theorem \ref{thm:main}. In subsection 3.2 we show the uniformly convergence. 
 \subsection{Long Time existence} Suppose $\Sigma_0=(x,u_0(x))$ is a graph of $u_0(x)$ in $N\PLH\R$ where $u_0(x)\in (u_0,u_1)$. Acoording to the standard theory of parabolic equation we can assume that $\Sigma_t$ is the smooth flow in \eqref{mct:flow} existing on $[0, T)$ for some $T>0$ and remaining as graphs over $N$ in $N\PLH\R$. We write $\Sigma_t$ as $(x,u(x,t))$ in $N\PLH\R$. Therefore according to Lemma \ref{lm:nonform} $u(x,t)$ satisfies \eqref{equ:non:form} with $u(x,0)=u_0(x)$ as follows.
  \be \label{equ:non:sform}
 \P_t u=g^{ij}u_{ij}+h(x,u)+g(x,u)\omega
 \ene 
  \bl\label{stepone}  With the assumptions in Theorem \ref{thm:main} $u(x,t)\in (u_0, u_1)$ for all $t\in [0, T)$.
  \el 
  \bp Otherwise there exists a time $t_0\in (0, T)$ and $x_0\in N$ such that $u(x,t)$ firstly achieves $u_0$ or $u_1$  at $p_0=(x_0, u(x_0,t_0))$. In the former case  we have 
   \be 
      Du(p_0)=0  \quad (\P_t u-g^{ij}u_{ij})(p_0)\leq 0
   \ene 
   However according to condition (1) and \eqref{equ:non:form} we always have 
   \be 
   (\P_t u-g^{ij}u_{ij})(p_0)=g(x_0,u_0)+h(x_0,u_0)\geq 0
   \ene 
   This is a contradiction to a strong maximum principle of parabolic equation. Therefore $u(x,t)$ never achieves $u_0$. With a similar derivation $u(x,t)$ never obtains $u_1$. 
   \ep 
 \bl\label{lm:gradient} There exists a positive constant $\mu_0$ independent of $T$ such that 
 \be 
  \omega (x,t) \leq \mu_0
 \ene
 for all $x\in N$ and $t\in [0, T)$. 
 \el 
 \bp Let $\eta =e^{Ku(x,t)}$ where $K$ is a positive constant determined later. We consider the evolution of $\eta\omega$ along the flow \eqref{equ:non:form}. According to \eqref{def:operator}, we have 
   \be \label{eq:step}
   L(\eta\omega)= (L\omega-\F{2}{\omega}g^{ij}\omega_i
   \omega_j)\eta+\omega L\eta 
   \ene 
   Now fix a positive constant $T_0$ where $T_0<T$. Assume $\eta\omega$ achieves the maximum on $N\times [0, T_0]$ at $(x_0, t_0)$. Furthermore we can assume $t_0>0$. Otherwise we can take larger $T_0$ until $T_0=T$. In this case nothing needs to prove. At this point $(x_0, t_0)$ combining \eqref{met} and \eqref{eq:step} one sees that 
    \be\label{mid:step}
   0\geq  \F{L\eta}{\eta}+|A|^2+Ric(v_N,v_N)+\psi_u\la v_N,v_N\ra+\F{u^k}{\omega^2}\psi_k
   \ene
   where $v_N=\F{Du}{\omega}$ and $\psi=-h(x,u)-g(x,u)\omega$.\\
   \indent  In the following we denote different constants independent of time $T_0$ by $C$. Since $N$ is closed and $\la v_N, v_N\ra \leq 1$, then $|Ric(v_N,v_N)|\leq C$. By condition (2) and Lemma \ref{stepone} we have 
               \begin{gather*}
               \psi_u\la v_N, v_N\ra=-\P_u(h(x,u))-\P_u(g(x,u))\omega \geq -\P_u(h(x,u))\geq C\\
           \F{u^k}{\omega^2}\psi_k= -\F{u^k}{\omega^2}(h_k(x,u)+g_k(x,u)\omega)\geq C
          \end{gather*}
          In summary at the point $(x_0, t_0)$ \eqref{mid:step} can be further simplified as follows. 
          \be \label{eq:std}
           0\geq  \F{L\eta}{\eta}-C
          \ene 
          where $C$ is independent of $T_0$. Now we expand $\F{L\eta}{\eta}$ and obtain that
  \be\label{eta:term}
    \F{L\eta}{\eta}= K^2g^{ij}u_iu_j+KLu
 \ene
   The first term in equation \eqref{eta:term} takes the form
   $$
   K^2g^{ij}u_iu_j=K^2(\sigma^{ij}-\F{u^i u^j}{1+|Du|^2})u_iu_j=K^2(1-\F{1}{\omega^2})
   $$
   Again by Lemma \ref{stepone} the second term in equation \eqref{eta:term} becomes
   \begin{align*}
     Lu &=\psi(x,u, Du)-\psi_{u_i}(x, u, Du)u_i)\\
        &=-h(x,u)-g(x,u)\F{1}{\omega}\geq C
   \end{align*}
  Putting these computations together into \eqref{eq:std} we obtain at point $(x_0,t_0)$ 
   \begin{align*}
   0\geq \F{L\eta}{\eta}-C&\geq K^2(1-\F{1}{\omega^2})-CK-C
   \end{align*}
   where $C$ are constants independent of $T_0$. Now choose $K_0$ sufficiently large only depending on $C$ we obtain that 
    \be 
      \omega (x_0, t_0)\leq \omega_0
    \ene
    Here $\omega_0>1$ is a positive constant only depends on $K_0$. Therefore $\omega_0$ is independent of $T_0$. Since $\eta\omega$ obtains its maximum on $N\times [0,t_0]$, then 
      \be 
      \omega \leq e^{K_0(u(x_0,t_0)-u(x,t))}\omega(x_0, t_0)\leq e^{K_0(u_1-u_0)}\omega_0
      \ene 
      where we apply Lemma \ref{lm:nonform}. Let $\mu_0$ be $K_0\omega_0$. Since $\mu_0$ does not depend on time $T_0$, we obtain the conclusion. 

\ep
  Lemma \ref{lm:nonform} and Lemma \ref{lm:gradient} give uniformly $C^0$ and $C^1$ bound of the flow in \eqref{equ:non:sform}. Therefore the flow \eqref{equ:non:sform} is uniformly parabolic. By the standard theory of nonlinear parabolic equation (see Section 5.5 in \cite{Kry87}), the flow \eqref{equ:non:sform} have long time existence. As a result the mean curvature type flow of graphs in \eqref{mct:flow} exists for all time and remains as graphs.
 \subsection{Uniformly Convergence} Now we show the uniformly convergence in Theorem \ref{thm:main}. From now on we assume $h(x,u)=h(u)$. We just show that $u(x,t)$ satisfies an uniformly parabolic equation \eqref{equ:non:sform} for all $t\in (0,\infty)$ with the uniformly $C^0$ and $C^1$ bound by Lemma  \ref{lm:nonform} and Lemma \ref{lm:gradient}. Thus we conclude that $|u_t(x,t)|\leq C$ where $C$ is independent of $C$. Moreover a key fact for uniformly convergence in our proof is given as follows.
  \bl \label{lm:integral}
   With the assumptions in Theorem \ref{thm:main}, it holds that
    \be
    \int_{0}^T\int_{M}|u_{t}|^2dxdt\leq C
    \ene for any $T>0$. Here $C$ is a positive constant independent of $T>0$.
 \el
 \bp Problem \eqref{equ:non:sform} is rewritten as
 \be
 u_t=div(\F{Du}{\omega})\omega+h(u)+g(x,u)\omega
 \ene
 where $div$ is the divergence fo $N$. Let $s(x,u)$ denote $e^{-\int h(u) du}$. 
 Applying the divergence formula one sees that
 \be\label{div:f1}
 \begin{split}
 	div(s(u)u_t\F{Du}{\omega})&=s(u)\la \F{Du}{\omega}, D(u_t)\ra+s'(u)u_t\F{|Du|^2}{\omega} + \\ &+s(u)(\F{u_t^2}{\omega}-g(x,u)u_t)-s(u)h(u)\F{u_t}{\omega};
 \end{split}
 \ene
 Then we compute the $t$-derivative of $s(u)\omega$ with equation \eqref{div:f1} as follows:
 \begin{align*}
 \P_t(s(u)\omega)&=s'(u)u_t\omega + s(u)\la \F{Du}{\omega}, Du_t\ra \\
 &=s'(u)\F{u_t}{\omega}+ div(s(u)u_t\F{Du}{\omega})-s(u)(\F{u_t^2}{\omega}-g(x,u)u_t)+s(u)h(u)\F{u_t}{\omega};\\
 &=div(s(u)u_t\F{Du}{\omega})-s(u)\F{u_t^2}{\omega}+s(u)g(x,u)u_t
 \end{align*}
 The last line follows since $s'(u)+s(u)h(u)\equiv 0$. Thus the $t$-derivative of $\int_{M}s(u)\omega dx$ takes the following form.
 \begin{align*}
 \P_t \int_{M}s(u)\omega dx=-\int_{M}s(u)\F{(u_t)^2}{\omega}dx+\P_t\int_{M}G(x,u)dx
 \end{align*}
 where $G(x,u)=\int s(u)g(x,u) du$. Reorganizing the above equation one sees that 
 \be
 \int_0^T\int_{M}s(u)\F{(u_t)^2}{\omega}dxdt=-
 \{\int_{M}s(u)\omega dx+\int_{M}G(x,u)dx\}|_0^T
 \ene
 Here $s(u)=e^{-\int h(u)du}$ and $G(x,u)=\int s(u)G(x,u)du$. The conclusion follows directly because from Lemma \ref{lm:nonform} and Lemma \ref{lm:gradient} $u(x,t)$ and $\omega$ are uniformly bounded. 
 \ep
Now we are ready to conclude the convergence of $u_t$ in Theorem \ref{thm:main}.
 \bl\label{lm:cn} Take  the assumptions in Theorem \ref{thm:main} and suppose $h(x,u)=h(u)$. Then $u_t(x,t)$ converges uniformly to 0.
 \el
 \bp Suppose $u_t(x,t)$ does not converge uniformly to 0. There are a constant $\Sc_0>0$, a sequence $\{x_n\}\in M$ and $\{t_{n}\}\rightarrow \infty$ such that
   $$
   |u_t(x_n, t_n)|\geq \Sc_0
   $$ Without loss of generality, we assume $t_{n+1}-t_{n}\geq 1$. Notice that in our settings all derivatives of $u(x,t)$ are uniformly bounded because of Lemma \ref{lm:nonform} and Lemma \ref{lm:gradient}. Thus after suitable choosing subsequence, there is a point $x_0\in M$ and $\Sc<\F{1}{2}$ such that
      $$ |u_t(x_0,t)|\geq \F{\Sc_0}{2} $$
      for all $t\in (t_n-\Sc, t_n+\Sc)$. This in turn implies that there is a compact domain $\Omega_{x_0}$ containing $x_0$ such that for all $x\in \Omega_{x_0}$ and $t\in (t_n-\Sc, t_n+\Sc)$, we have
           $$ |u_t(x,t)|\geq \F{\Sc_0}{4} $$
      because all higher derivatives of $u(x,t)$ are uniformly bounded.
  Thus
       $$
       \int_{\Omega_{x_0}}\int_{t_n-\Sc}^{t_{n}+\Sc}|u_t|^2 dtdx\geq \F{\Sc_0}{2}vol(\Omega_{x_0})2\Sc>0
       $$
  as $t_n\rightarrow \infty$. This gives a contradiction to Lemma \ref{lm:integral}. Hence $u_t(x,t)$ converges uniformly to 0. \\
  \ep
  Now we conclude the uniformly convergence of Theorem \ref{thm:main}. Since  $u(x,t)$ and all derivatives of $u(x,t)$ are uniformly bounded. Then there exists a sequence $\{t_i\}_{i=1}^\infty$ going to infinity such that $u(x, t_i)$ converges to a smooth function $u_{\infty}(x)$. Because $u(x,t)$ satisfies that 
  $$
    u_t = g^{ij}u_{ij}+h(u)+g(x,u)\omega =-\omega (H(u)-h(u)\la \vec{v},\P_r\ra+g(x,u))
  $$
   and by Lemma \ref{thm:main} $u_t\rightarrow 0$ as $t\rightarrow \infty$, $u_\infty(x)$ is the solution to the problem 
     \be 
      H(u)= h(u)\la \vec{v},\P_r\ra-g(x,u)
     \ene 
 The proof of Theorem \ref{thm:main} is complete.

\section{The proof of Corollary \ref{cor:one} and Corollary \ref{cor:two}}\label{section:G}
In this section we apply Theorem \ref{thm:main} into the setting of warped product manifolds. Those applications are based on Lemma \ref{lm:transform} which gives a one-to-one correspondence between graphs in warped product manifolds and those in product manifolds. 
\subsection{The proof of Corollary \ref{cor:one}} For convenience we restate corollary \ref{cor:one} as follows. 
  \begin{cor}(Corollary \ref{cor:one}) Suppose $N$ is a n-dimensional closed Riemannian manifold. Let $\phi(u)$ be a positive smooth function on $\R$ and $f(x,u)$ be a smooth function where $x\in N$ and $u\in \R$. Suppose the following two conditions hold:
 	\begin{enumerate}
 		\item  for all $x\in N$ there exists two constants $u_0<u_1$ satisfying 
 		\be\label{con:bars}
 		f(x,u_0)\geq n\F{\phi'(u_0)}{\phi(u_0)}\quad  f(x,u_1)\leq n\F{\phi'(u_1)}{\phi(u_1)}
 		\ene
 		\item moreover $\P_u (f(x,u)\phi(u))\leq 0$ for all $x\in N$ and $u\in [u_0,u_1]$. 
 	\end{enumerate}
 	Then there exists a  smooth hypersurface $\Sigma=(x,u(x))$ in warped product manifold $N\PLH_{\phi}\R$ such that its mean curvature is $f(x,u(x))$. 
 \end{cor}

\bp Now we work in the product manifold $N\PLH\R$. Let $g(x,\vp)=f(x,u)\phi(u)$ and $h(\vp)=-\phi'(u)$ where $u(x)=\Phi^{-1}(\vp(x))$ from \eqref{relation}. As a result $\F{\P\vp}{\P u}=\F{1}{\phi(r)}>0$. It is easy to see that condition (1) and condition (2) imply the two conditions in Theorem \ref{thm:main}. Theorefore Theorem \ref{thm:main} implies that there is a smooth function $\vp(x)$ such the mean curvature of its graph in $N\PLH\R$ is 
   \begin{align*}
     H(\vp) &=-\phi'(u)\la\vec{v},\P_r\ra +f(x,u)\phi(u)\\
            &=\phi'(u)\F{1}{\omega}+f(x,u)\phi(u)
   \end{align*} 
 where $u=\Phi^{-1}(\vp(x))$ and $\omega=\sqrt{1+|D\vp|^2}$. Lemma \ref{lm:transform} we obtain the conclusion. 
\ep
\subsection{The proof of Corollary \ref{cor:two}} We restate Corollary \ref{cor:two} as follows. 
\begin{cor} ( Corollary \ref{cor:two}) Let $a<b$ be two fixed constants. Suppose a warped function $\phi(u)$ satisfies
	\be\label{eq:ust}
	\phi'(a)<0\quad \phi'(b)>0 \quad \phi''(u)\geq 0\quad \text{for all $u\in (a,b)$}
	\ene
	Let $N$ be a closed Riemananian manifold. Suppose $\Sigma_0$ is a graphical smooth hypersurface in a warped product manifold $N\PLH_{\phi}(a,b)$. Consider a weighted mean curvature flow 
	\be \label{steu}
	\P_t F(p)=-\phi^2(u)H\vec{v}
	\ene 
	where $p=(x,u)\in N\PLH_{\phi}\R$, $F(0,p)=p$ for $p\in \Sigma$ and $H$ is the mean curvature of $F_t(\Sigma_0)$. Then the flow of $\Sigma_0$ in \eqref{steu} exists for all times, remains graphical and converges uniformly to the totally geodesic slice in $N\PLH_{\phi}(a,b)$. 
\end{cor}
\bp Suppose $\Sigma_t$ is a family of graphs satisfying the flow in \eqref{steu} with the initial graph $\Sigma_0$, represented as $(x,r(x,t))$. Then according to Lemma \ref{lm:ef} $r(x,t)$ satisfies 
\be\label{stded}
\P_t r(x,t)=\phi(r)\F{1}{\omega}(g^{ij}\vp_{ij}-n\phi'(r))
\ene
where $\vp$ is the modified height satisfying $\vp'(r)=\F{1}{\phi(r)}$ (see \eqref{relation}). Then \eqref{stded} is equivalent to 
\be\label{eust}
\P_t \vp=\F{1}{\omega}g^{ij}\vp_{ij}-n\F{\phi'(r(\vp))}{\omega}
\ene 
where $r(\vp)$ is the inverse function of $\vp=\vp(r)$ (see \eqref{relation}). Let $\vp_0=\vp(a)$ and $\vp_1=\vp(b)$. Therefore $\vp_1>\vp_0$ by the definition of $\vp(r)$. 
Consider the flow $(x,\vp(x,t))$ in the product manifold $N\PLH [\vp_0,\vp_1]$. Let $g(x,\vp)=0$ and $h(x,\vp)=-n\phi'(r(\vp))$. From \eqref{eq:ust}, we have 
  \begin{gather} 
  h(x,\vp_0)+g(x,\vp_0)=-n\phi'(a)\geq 0\quad \\ h(x,\vp_1)+g(x,\vp_1)=-n\phi'(b)\leq 0
  \end{gather}
  The above verify the two conditions in Theorem \ref{thm:main}.
 Together with Lemma \ref{lm:nonform} and Theorem \ref{thm:main}, we obtain the long time existence of the flow \eqref{steu} for a smooth graph $\Sigma_0$ in $N\PLH_{\phi}[a,b]$. \\
 \indent Now we show the uniformly convergence. 
 First notice that by condition \eqref{eq:ust} there is only one point $x_0\in (a,b)$ such that $\phi'(x_0)=0$. By the proof of Lemma \ref{lm:ef} the slice $N\PLH\{x_0\}$ is totally geodesic in $N\PLH_{\phi}[a,b]$.\\
 \indent According to Lemma \ref{lm:gradient}  and Lemma \ref{stepone} $\vp(x,t)$ has uniformly $C^0$ and $C^1$ bound. All derivatives of $\vp(x,t)$ are uniformly bounded and $\vp_t$ converges to $0$ uniformly. If we can show that $\vp(x,t)$ converges to $x_0$ uniformly, the conclusion follows. \\
 \indent Now we argue that the above claim is true.  Next we consider the flow in \eqref{steu} of each slice $N\PLH\{r_0\}$ where $r_0\in [a,b]$. It should satisfy an ODE equation as follows:
 \be 
 \P_t r(x,t)=-n \phi'(r(x,t))
 \ene 
 where $r(x,0)=r_0\in (a,b)$. From condition \eqref{eq:ust} the flow in \eqref{steu} of each slice is a family of slices and will converge uniformly to the slice $N\PLH\{x_0\}$. One the other hand, from Huisken \cite{Hui86} the flow in \eqref{stded} has a disjoint principle. That is if the initial closed hypersurfaces are disjoint with each other, their flows are always disjoint with each other at any time as long as they exist smoothly. Thus the flow of $\Sigma_0$ is always disjoint with the flows of $N\PLH\{a\}$ and $N\PLH\{b\}$. Now all three flows exist for all time and the latter two flows converge to $N\PLH \{x_0\}$ uniformly. Therefore so is the flow of $\Sigma_0$. We obtain the claim at the beginning.\\
 \indent The proof of Corollary \ref{cor:two} is complete. 
\ep
 
\section*{Acknowledgment}
The first author is supported by the National Natural Science Foundation of China (Grant No. 11401578). The second author is very appreciated to the encouragement from Professor Lixin Liu and Professor Zheng Huang. The second author is supported by the National Natural Science Foundation of China (Grant No. 11261378 and Grant No. 11521101).

\bibliographystyle{abbrv}	
\bibliography{Ref_Thesis}
\end{document}